\documentstyle[12pt]{article}
\title{On the babylonian method of extracting root squares}
\author{Vilma A. S. Sant'Anna\thanks{Col\'egio Martinus J\'unior, Rua Arthur Mohr, 144, 81070-040 Curitiba, PR, Brazil} \and Adonai S. Sant'Anna\thanks{To whom correspondence should be sent: C.P. 019081, Dep. Matem\'atica, Universidade Federal do Paran\'a, Curitiba, PR, 81531-990, Brazil}}
\date{ }
\begin{document}

\newtheorem{definicao}{\bf Definition}
\newtheorem{teorema}{\bf Theorem}
\newtheorem{lema}{\bf Lemma}
\setlength{\unitlength}{1mm}
\renewcommand{\thefootnote}{\arabic{footnote}}
\setcounter{footnote}{0}
\newcounter{cms}
\setlength{\unitlength}{1mm}

\maketitle

\begin{abstract}
We discuss the babylonian method of extracting the root square of a number, from the point of view of modern mathematics. We also speculate that the babylonian mathematics was rich enough for a generalization of this method, despite the lack of general statements and justified procedures in their mathematics.
\end{abstract}

\noindent
Keywords: Babylonian mathematics; Newton's method; root squares.\\\\
AMS Subject Classification: 01A17

\section{Introduction}

	The oldest babylonian mathematical texts known to us date from the period 1900-1600 b.C.. It is well known, for example, that this people knew how to extract the root square of any positive number. Despite the fact that there was no general statements, rules or justified procedures in their mathematics, it seems reasonable to suppose that their probable line of thought was the following \cite{Bunt-88}:

\begin{enumerate}

\item The root square of, say, 17 is a number whose square is 17.

\item The root square of 17 is approximatelly (as a first approximation) 4. Let us call this number $r_1$.

\item The number $r_1 = 4$ is not the root square of 17; but if we multiply it by $\frac{17}{4}$, we have 17 as result. In other words, both numbers, 4 and $\frac{17}{4}$ are good approximations of the root square of 17.

\item The arithmetic mean of 4 and $\frac{17}{4}$ should be a better approximation for the root square of 17. This arithmetic mean, which is a second approximation, is $4\frac{1}{8}$. We refer to this number as $r_2$.

\item The number $r_2 = 4\frac{1}{8}$ is not the root square of 17; but if we multiply it by $\frac{17}{r_2}$, we get 17 as result. In other words, both numbers, $r_2$ and $\frac{17}{r_2}$, are good approximations of the root square of 17.

\item The arithmetic mean of $r_2$ and $\frac{17}{r_2}$ should be $r_3$, and so on.

\end{enumerate}

	It seems quite obvious that this procedure is a rudiment of a numerical method. If we consider the case where the number of approximations goes to infinity, we have the limit $L$, which is the square root of 17, satisfying the following:

\begin{equation}
\frac{L+\frac{17}{L}}{2} = L,
\end{equation}

\noindent
which gives the equation $L^2 = 17$.

\section{Newton's Method}

	Now, consider the Newton's method of approximating to the zeroes of a given function $f(x)$:

\begin{equation}
x_{n+1} = x_n-\frac{f(x_n)}{f'(x_n)},\label{Newton}
\end{equation}

\noindent
where $f'(x)$ is the first derivative of $f(x)$ with respect to $x$, and $x_1$ is a first approximation to a root of the equation $y = f(x)$.

	In the particular case

\begin{equation}
f(x) = x^2-17,
\end{equation}

\noindent
where one of the roots is $\sqrt{17}$, we have the following:

\begin{equation}
x_{n+1} = \frac{x_n+\frac{17}{x_n}}{2},\label{quadrado}
\end{equation}

\noindent
which is exactly the same method used by the babylonians. In other words, the babylonians used a particular case of Newton's method although they were not aware about that.

\section{Speculating}

	Now, let us make some speculation. Was it possible, for a babylonian mathematician, to use very simple arguments (similar to those presented in the first Section of this paper) in order to extract the root cubic of a positive number? We believe that the answer is positive.

	The root cubic of 17 is one of the zeroes of

\begin{equation}
f(x) = x^3 - 17.
\end{equation}

	According to Newton's method given by equation (\ref{Newton}) the approximations are obtained as it follows:

\begin{equation}
x_{n+1} = \frac{x_n+x_n+\frac{17}{x_n^2}}{3},\label{cubo}
\end{equation}

\noindent
where $x_1$ is the first approximation. If we were babylonian mathematicians we could make some analogy with the method presented in the first Section as it follows:

\begin{enumerate}

\item The root cubic of 17 is a number whose cubic is 17.

\item The root cubic of 17 is approximatelly (as a first approximation) 2. Let us call this number $r_1$.

\item The number $r_1 = 2$ is not the root cubic of 17; but if we perform the product $2\times 2 \times \frac{17}{2^2}$, we have 17 as a result. In other words, the numbers, 2 and $\frac{17}{2^2}$ are good approximations of the root cubic of 17.

\item The arithmetic mean of 2, 2, and $\frac{17}{2^2}$ should be a better approximation for the root cubic of 17. This arithmetic mean, which is a second approximation, is $2\frac{3}{4}$. We call this number as $r_2$.

\item The number $r_2 = 2\frac{3}{4}$ is not the root cubic of 17; but if we calculate the product $r_2\times r_2\times \frac{17}{r_2^2}$, we get 17 as result. In other words, the numbers, $r_2$ and $\frac{17}{r_2^2}$, are good approximations of the root cubic of 17.

\item The arithmetic mean of $r_2$, $r_2$, and $\frac{17}{r_2^2}$ should be $r_3$, and so on.

\end{enumerate}

	In other words, we are again using the notion of arithmetic mean, just as in the case of root squares. So, it seems natural that a hypothetical clever babylonian mathematician could create a method for solving the problem of the root cubic of a positive number, although there was (obviously) no knowledge about Newton's method at that time.

\section{A Simple Generalization}

	The reader can easily prove, that a generalization of equations (\ref{cubo}) and (\ref{quadrado}) is given by:

\begin{equation}
x_{n+1} = \frac{\sum_{i=1}^{m-1}x_n+\frac{r}{x_n^{m-1}}}{m},
\end{equation}

\noindent
where $x_1$ is the first approximation of $\sqrt[m]{r}$ and $\sum_{i=1}^{m-1}x_n = (m-1)x_n$. It can be easily proved as well that any positive real number can be the first approximation $x_1$.

\end{document}